\documentclass[a4paper,UTF8]{article}
\usepackage[a4paper,left=3cm,right=3cm,top=3cm,bottom=3cm]{geometry}
\usepackage{amsmath} 
\usepackage{amssymb} 
\usepackage{amsfonts} 
\usepackage{graphicx} 
\usepackage{xcolor} 
\usepackage{BOONDOX-cal}
\usepackage{amsthm}
\usepackage{tikz}
\usepackage{authblk}
\usepackage{fancyhdr}
\theoremstyle{plain}
\newtheorem{theorem}{Theorem}[section]
\newtheorem*{theorem*}{Theorem}
\newtheorem{proposition}[theorem]{Proposition}
\newtheorem*{proposition*}{Proposition}

\newtheorem*{lemma*}{Lemma}

\newtheorem*{corollary*}{Corollary}
\newtheorem{definition}[theorem]{Definition}
\newtheorem*{definition*}{Definition}

\newtheorem*{remark*}{Remark}
\newtheorem{example}[theorem]{Example}
\newtheorem*{example*}{Example}

\usepackage{hyperref}
\usepackage[numbers,sort&compress]{natbib}

\begin{document}
\title{The generalized upper box dimension}

\author[1]{Lipeng Wang \thanks{Corresponding author. \\ E-mail address: lipengwang@jou.edu.cn, wxli@math.ecnu.edu.cn}}
\author[2]{Wenxia Li}

\affil[1]{School of Science, Jiangsu Ocean University, 59 Cangwu Rd., Lianyungang 222005, P. R. China}

\affil[2]{School of Mathematical Sciences, East China Normal University, 500 Dongchuan Rd., Shanghai 200241, P. R. China}

\renewcommand*{\Affilfont}{\small\it}
\date{}
\maketitle

\noindent\textbf{Abstract.}
We introduce the generalized upper box dimension which is defined for any set, whether the set is bounded or unbounded. We study basic properties of the generalized upper box dimension. We prove that the generalized upper box and upper box dimensions coincide for bounded sets. Furthermore, we also show that the modified generalized upper box dimension equals the packing dimension. So the generalized upper box dimension can be seen as a reasonable generalization of the upper box dimension. As an application, we prove the generalized upper box dimension is zero if and only if the quasi-Assouad dimension is zero. We also show that the upper spectrum is of full dimension is equivalent to the Assouad spectrum is of full dimension and the upper spectrum is zero is equivalent to the Assouad spectrum is zero.



\section{Introduction}
The upper box dimension is an important fractal dimension in fractal geometry. Compared to many fractal dimensions such as the Hausdorff, packing and Assouad dimensions, the upper box dimension is special because the upper box dimension is only defined for bounded sets and not for unbounded sets. When studying problems concerning the upper box dimension, we need to assume that the set being studied is bounded, which is inconvenient. In \cite{FHHTY}, Fraser et al. said that "When we discuss the upper box-counting dimension we are implicitly referring to bounded sets only, since the definition does not readily apply to unbounded sets". A natural problem is how to define the generalized upper box dimension. This generalized upper box dimension is defined for any set, whether the set is bounded or unbounded. For bounded sets, the generalized upper box dimension can have the basic properties of the upper box dimension. We try to use the Assouad and upper spectra to define the generalized upper box dimension.

The Assouad and upper spectra are both variations of the Assouad dimension. When Assouad studied embedding theory \cite{A1, A2, A3}, the Assouad dimension was popularised. The Assouad dimension has attracted great interest in fractal geometry \cite{BKR, Fr1, Fr2, FHOR, LWX, LLMX, O, YD}.
In the past few years, research on variations of the Assouad dimension has become a focus in fractal geometry \cite{CDW, CWC, FHHTY, Fr1, Fr2, FY2, LX, XDW}. L\"{u} and Xi \cite{LX} introduced the quasi-Assouad dimension which is a quasi-Lipschitz invariant. Xi et al. \cite{XDW} studied  Assouad-minimality of a class of Moran sets under quasi-lipschitz
mappings and used the quasi-Assouad dimension to study Assouad-minimality of another class of Moran sets.
 Fraser and Yu \cite{FY2} introduced the Assouad spectrum which interpolates between the upper box and quasi-Assouad dimensions. Fraser et al. \cite{FHHTY} studied the relationship between the upper spectrum and the Assouad spectrum.

In this article, we use the Assouad and upper spectra to introduce the generalized upper box dimension which can be seen as a reasonable generalization of the upper box dimension. We study basic properties and application of the generalized upper box dimension.

\section{Preliminaries}
Let $E\subset\mathbb{R}^d$ be nonempty and bounded. For any $\delta>0$, the smallest number of closed balls of radius $\delta>0$ required to cover $E$ is denoted by $N_\delta(E)$. We use $\parallel\cdot\parallel$ to denote the Euclidean norm in $\mathbb{R}^d$ and use $B(x, R):=\{y\in\mathbb{R}^d:\parallel y-x\parallel\leq R\}$ to denote the closed ball with a centre $x\in\mathbb{R}^d$ and radius $R>0$.

\subsection{Upper box dimension}
\begin{definition}[\cite{F}]
For any bounded $F\subset\mathbb{R}^d$, the upper box dimension of $F$, denoted by $\overline{\dim}_{B} F$, is defined by
$$\overline{\dim}_{B} F:=\limsup\limits_{\delta\rightarrow 0}\frac{\log N_\delta(F)}{-\log\delta},$$
where $\overline{\dim}_{B} \emptyset:=0.$
\end{definition}

\subsection{Packing and modified upper box dimensions}
The following can be found in Falconer's textbook \cite{F}.

\noindent For any $F\subset\mathbb{R}^d$, $s\geq 0$ and $\delta>0$, define
\begin{align*}
  \mathcal{P}_\delta^s :=\sup\bigg\{\sum\limits_{i=1}^{+\infty}|B_i|^s:&\ \{B_i\}\ \mbox{is a at most countable collection of disjoint closed balls} \\
   & \mbox{of radii at most}\ \delta\ \mbox{with centres in}\ F\bigg\}.
\end{align*}
Since $\mathcal{P}_\delta^s$ decreases as $\delta$ decreases to $0$, the limit
$$\mathcal{P}_0^s:=\lim\limits_{\delta\rightarrow 0}\mathcal{P}_\delta^s$$ exists, possibly $0$ or $+\infty$.
Next define an outer measure
$$\mathcal{P}^s(F):=\inf\bigg\{\sum\limits_{i=1}^{+\infty}\mathcal{P}_0^s(F_i): F\subset\bigcup\limits_{i=1}^{+\infty} F_i\bigg\}.$$
\begin{definition}[\cite{F}]
The packing dimension of $F$, denoted by $\dim_P F$, is defined by
\begin{align*}
  \dim_P F & :=\sup\{s\geq 0: \mathcal{P}^s(F)=+\infty\} \\
   & =\inf\{s\geq 0: \mathcal{P}^s(F)=0\},
\end{align*}
where $\dim_P \emptyset:=0.$
\end{definition}
\begin{definition}[\cite{F}]
The modified upper box dimension of $F$, denoted by $\overline{\dim}_{MB} F$, is defined by
\begin{align*}
  \overline{\dim}_{MB} F & :=\inf\bigg\{\sup\limits_{i\geq 1}\overline{\dim}_B F_i: F\subset\bigcup\limits_{i=1}^{+\infty}F_i,\ \mbox{where}\ F_i\ \mbox{is bounded for each}\ i\in \mathbb{N}\bigg\}\\
   & =\inf\bigg\{\sup\limits_{i\geq 1}\overline{\dim}_B F_i: F=\bigcup\limits_{i=1}^{+\infty}F_i,\ \mbox{where}\ F_i\ \mbox{is bounded for each}\ i\in \mathbb{N}\bigg\}.
\end{align*}
\end{definition}

\begin{proposition}[\cite{F} Proposition 3.9] \label{pdmubd}
For $F\subset \mathbb{R}^d$, we have
  $$\dim_P F=\overline{\dim}_{MB} F.$$
\end{proposition}

\subsection{Assouad and upper spectra and Assouad and quasi-Assouad dimensions}

\begin{definition} [\cite{Fr1}]
For $F\subset\mathbb{R}^d$, the Assouad dimension of $F$, denoted by $\dim_A F$, is defined by
\begin{align*}
\dim_A F:=\inf\bigg\{s\geq 0:\ & \mbox{there exist constants} \ C>0,\ \rho>0\ \mbox{such that, for all}\\
                      & 0<r<R<\rho\ \mbox{and}\ x\in F,\ N_r(B(x,R)\cap F)\leq C\left(\frac{R}{r}\right)^s\bigg\},
\end{align*}
where $\dim_A \emptyset:=0.$
\end{definition}

\begin{definition}[\cite{FHHTY, LX}] For any $\theta\in (0, 1)$ and $F\subset\mathbb{R}^d$, the upper spectrum of $F$, denoted by $\overline{\dim}_{A}^\theta F$, is defined by
\begin{align*}
\overline{\dim}_{A}^\theta F :=\inf\bigg\{s\geq 0:\ & \mbox{there exist a constant} \ C>0\ \mbox{such that, for all}\\
                                           & 0<r\leq R^{1/\theta}<R<1\ \mbox{and}\ x\in F,\ N_r(B(x,R)\cap F)\leq C\left(\frac{R}{r}\right)^s\bigg\}.
\end{align*}

\noindent The quasi-Assouad dimension of $F$, denoted by $\dim_{qA} F$, is defined by $$\dim_{qA} F:=\lim\limits_{\theta\rightarrow 1}\overline{\dim}_{A}^\theta F,$$
where $\overline{\dim}_{A}^\theta \emptyset:=0$ for each $\theta\in (0, 1)$.

\end{definition}

\begin{definition} [\cite{FY2}] For any $\theta\in (0, 1)$ and $F\subset\mathbb{R}^d$, the Assouad spectrum of $F$, denoted by $\dim_{A}^\theta F$, is defined by
\begin{align*}
\dim_{A}^\theta F :=\inf\bigg\{s\geq 0:\ & \mbox{there exist a constant} \ C>0\ \mbox{such that, for all}\\
                                & 0<r=R^{1/\theta}<R<1\ \mbox{and}\ x\in F,\ N_r(B(x,R)\cap F)\leq C\left(\frac{R}{r}\right)^s\bigg\},
\end{align*}
where $\dim_{A}^\theta \emptyset:=0$ for each $\theta\in (0, 1)$.
\end{definition}

\noindent The following propositions are important properties of the Assouad and upper spectra. We will use them later to prove some of our results.

\begin{proposition}[\cite{FY2}]\label{asbp}
(1) Let $E\subset F\subset\mathbb{R}^d$ and $\theta\in (0,1)$. Then $\dim_{A}^\theta E\leq \dim_{A}^\theta F.$

\noindent (2) For any $E, F\subset\mathbb{R}^d$ and $\theta\in (0,1)$, we have $\dim_{A}^\theta E\cup F=\max\{\dim_{A}^\theta E, \dim_{A}^\theta F\}.$

\noindent (3) Let $E, F\subset\mathbb{R}^d$ and $\theta\in (0,1)$. If there exists a bi-Lipschitz mapping $f$ such that $f(E)=F$, then $\dim_{A}^\theta E= \dim_{A}^\theta F.$

\noindent (4) Let $F\subset\mathbb{R}^d$ and $\theta\in (0,1)$. Then $\dim_{A}^\theta F= \dim_{A}^\theta \overline{F},$ where $\overline{F}$ denote the closure of $F$ in $\mathbb{R}^d$.

\noindent (5) For any $E\subset\mathbb{R}^m, F\subset\mathbb{R}^n$ and $\theta\in (0,1)$, we have
$\dim_{A}^\theta E\times F\leq \dim_{A}^\theta E+\dim_{A}^\theta F.$
\end{proposition}

\begin{proposition} [\cite{FHHTY}] \label{usasbp}

\noindent (1) For any $F\subset\mathbb{R}^d$, $\overline{\dim}_{A}^\theta F$ is nondecreasing in $\theta$.

\noindent (2) For any $F\subset\mathbb{R}^d$, $\dim_{A}^\theta F\leq \overline{\dim}_{A}^\theta F\leq \dim_{qA} F\leq \dim_A F$ for each $\theta\in (0,1)$.
\end{proposition}

\begin{proposition}[\cite{Fr1} Lemma 3.4.4, \cite{FY2} Corollary 3.2]\label{asub}
For any bounded $F\subset\mathbb{R}^d$,
$$\overline{\dim}_{B} F \leq \dim_{A}^\theta F\leq \min\left\{\frac{\overline{\dim}_B F}{1-\theta}, \dim_{qA} F\right\}$$
and
$$\lim\limits_{\theta\rightarrow 0}\dim_{A}^\theta F=\overline{\dim}_B F.$$
\end{proposition}

\begin{proposition}[\cite{FHHTY} Theorem 2.1]\label{usas}
For any $\theta\in (0, 1)$ and $F\subset\mathbb{R}^d$,
$$\overline{\dim}_{A}^\theta F=\sup\limits_{0<\theta^{'}\leq \theta}\dim_{A}^{\theta^{'}} F.$$
\end{proposition}

\begin{proposition}[\cite{FY2} Corollary 3.5]\label{asc}
For any $F\subset\mathbb{R}^d$, the function $\theta\mapsto\dim_{A}^\theta F$ is continuous in $\theta\in (0, 1)$.
\end{proposition}

\begin{proposition}[\cite{FY2} Corollary 3.6]\label{asa}
For any $F\subset\mathbb{R}^d$, if for some $\theta\in (0, 1)$, we have $\dim_{A}^\theta F=\dim_A F$, then
$$\dim_{A}^{\theta^{'}} F=\dim_A F$$
for all $\theta^{'}\in [\theta, 1)$.
\end{proposition}

\begin{proposition}[\cite{Fr1} Theorem 3.3.1]\label{ai}
For any $F\subset\mathbb{R}^d$ and $0<\theta_1<\theta_2<1$, we have
\begin{align*}
  \left(\frac{1-\theta_2}{1-\theta_1}\right)\dim_{A}^{\theta_2} F & \leq \dim_{A}^{\theta_1} F \\
   & \leq  \left(\frac{1-\theta_2}{1-\theta_1}\right)\dim_{A}^{\theta_2} F+ \left(\frac{\theta_2-\theta_1}{1-\theta_1}\right)\dim_{A}^{\theta_1/\theta_2} F.
\end{align*}

\end{proposition}

\section{Results}

\begin{proposition}\label{asusl}
For any $F\subset\mathbb{R}^d$, we have
$$\limsup_{\theta\rightarrow 0}\dim_{A}^\theta F=\lim_{\theta\rightarrow 0}\overline{\dim}_{A}^\theta F.$$
\end{proposition}

\begin{definition}\label{gubdd}
For any $F\subset\mathbb{R}^d$, the generalized upper box dimension of $F$, denoted by $\overline{\dim}_{GB} F$, is defined by
\begin{align*}
  \overline{\dim}_{GB} F :&=\limsup_{\theta\rightarrow 0}\dim_{A}^\theta F \\
   & =\lim_{\theta\rightarrow 0}\overline{\dim}_{A}^\theta F.
\end{align*}
\end{definition}

\begin{proposition}\label{gubdubd}
Let $F\subset\mathbb{R}^d$ be bounded. Then
$$\overline{\dim}_{GB} F=\overline{\dim}_B F.$$
\end{proposition}

\noindent For the unbounded set $F\subset\mathbb{R}^d$, its upper box dimension can usually be defined as follows:
$$\overline{\dim}_{GB}^* F:=\lim_{R\rightarrow +\infty}\overline{\dim}_B(B(O, R)\cap F),$$
where $B(O, R)$ denote the closed ball with a center at the origin $O$ in $\mathbb{R}^d$ and radius $R>0$.

\noindent It is easy to see that
$$\overline{\dim}_{GB}^* F=\lim_{R\rightarrow +\infty}\overline{\dim}_B(B(O, R)\cap F)=\lim_{R\rightarrow +\infty}\overline{\dim}_{GB}(B(O, R)\cap F)\leq \overline{\dim}_{GB} F.$$

\noindent Next Example \ref{gubr} shows that there exists a set $E\subset\mathbb{R}$ such that $\overline{\dim}_{GB}^* E<\overline{\dim}_{GB} E.$

\noindent By Propositions \ref{usasbp} and \ref{asub}, we have
$\overline{\dim}_{B} F\leq \dim_{A}^\theta F\leq \overline{\dim}_{A}^\theta F\leq \dim_{qA} F\leq \dim_A F$ for each bounded set $F\subset\mathbb{R}^d$ and $\theta\in (0,1)$.
As $\dim_{qA} F:=\lim\limits_{\theta\rightarrow 1}\overline{\dim}_{A}^\theta F$,
we believe that the generalized upper box dimension should be defined as $\limsup\limits_{\theta\rightarrow 0}\dim_{A}^\theta F$ to make it more consistent with other fractal dimensions. Next Proposition \ref{gubdas} shows that $\lim\limits_{\theta\rightarrow 0}\dim_{A}^\theta F$ exists.

\begin{example}\label{gubr}
For any positive integer $n\geq 2$, let $\delta_n:=\frac{1}{n^{1+1/(n-1)}}$ and $E:=\bigcup\limits_{n=2}^{+\infty}\{n+i\delta_n:i=0,1,\cdots, n\}.$
Then $$\overline{\dim}_{GB} E=1>0=\overline{\dim}_{GB}^* E.$$
\end{example}

In fractal geometry, the fractal dimension can reflect the complexity of a set. Some fractal dimension of a set is zero,  we believe that the set is not complex. In Example \ref{gubr}, $\overline{\dim}_{GB}^* E=0$, then we have the set $E$ is not complex. In fact, as it approaches infinity, the set $E$ contains increasingly dense arithmetic sequences, so when viewed as a whole, the set is complex. The $\overline{\dim}_{GB}^*$ cannot distinguish this complexity, so from the perspective of distinguishing complexity, the $\overline{\dim}_{GB}$ is more refined than the $\overline{\dim}_{GB}^*$.

\begin{theorem}\label{gubdasus}
For any $\theta\in (0, 1)$ and $F\subset\mathbb{R}^d$, we have
$$\overline{\dim}_{GB} F \leq \dim_{A}^\theta F\leq \min\left\{\frac{\overline{\dim}_{GB} F}{1-\theta}, \dim_{qA}F\right\}$$
and
$$\overline{\dim}_{GB} F \leq \overline{\dim}_{A}^\theta F\leq \min\left\{\frac{\overline{\dim}_{GB} F}{1-\theta}, \dim_{qA}F\right\}$$
\end{theorem}

\begin{proposition}\label{gubdas}
For any $F\subset\mathbb{R}^d$, we have
$$\overline{\dim}_{GB} F =\lim_{\theta\rightarrow 0}\dim_{A}^\theta F.$$
\end{proposition}

\begin{proposition}\label{gubdp}
(1) Let $E\subset F\subset\mathbb{R}^d$. Then $\overline{\dim}_{GB} E\leq \overline{\dim}_{GB} F.$

\noindent (2) For $E, F\subset\mathbb{R}^d$, we have $\overline{\dim}_{GB} E\cup F=\max\{\overline{\dim}_{GB} E, \overline{\dim}_{GB} F\}.$

\noindent (3) Let $E, F\subset\mathbb{R}^d$. If there exists a bi-Lipschitz mapping $f$ such that $f(E)=F$, then $\overline{\dim}_{GB} E= \overline{\dim}_{GB} F.$ Recall that the mapping $f:E\rightarrow\mathbb{R}^d$ is called bi-Lipschitz if there exists $L\geq 1$ such that
$$L^{-1}\|x-y\|\leq \|f(x)-f(y)\|\leq L\|x-y\|,$$
for each $x, y \in E$.

\noindent (4) Let $F\subset\mathbb{R}^d$. Then $\overline{\dim}_{GB} F= \overline{\dim}_{GB} \overline{F},$ where $\overline{F}$ denote the closure of $F$ in $\mathbb{R}^d$.

\noindent (5) For any $E\subset\mathbb{R}^m, F\subset\mathbb{R}^n$, we have
$\overline{\dim}_{GB} E\times F\leq \overline{\dim}_{GB} E+\overline{\dim}_{GB} F.$

\noindent (6) For any $F\subset\mathbb{R}^d$, we have $\dim_P F\leq \overline{\dim}_{GB} F\leq \dim_A F.$
\end{proposition}

As an application, we will use properties of the generalized upper box dimension to prove the following Theorems \ref{gubdqad} and \ref{ua}.

\begin{theorem}\label{gubdqad}
Let $F\subset\mathbb{R}^d$. Then
\noindent $$\overline{\dim}_{GB} F=0\ \mbox{if and only if}\ \dim_{qA} F=0.$$
\end{theorem}

\noindent\textbf{Remark.} In Theorem 2.7 in \cite{WL}, for any $0<b<c\leq 2$, we have proved that there exists a set defined by digit restrictions $E_{S, D}\subset\mathbb{R}^2$ which is bounded such that $\overline{\dim}_{GB} E_{S, D}=\overline{\dim}_B E_{S, D}=b<\dim_{qA} E_{S, D}=c.$ So
$\dim_{qA} F=c\in (0, d]$ cannot imply $\overline{\dim}_{GB} F=c$ and vice versa. When $F$ is bounded, the equivalence in Theorem \ref{gubdqad} is known in Corollary 2.3 in \cite{FHHTY}.

\begin{theorem}\label{ua}
Let $F\subset\mathbb{R}^d$. Then, for any $\theta\in (0,1)$, we have
$$\overline{\dim}_{A}^\theta F=d\ \mbox{if and only if}\ \dim_{A}^\theta F=d,$$
and
$$\overline{\dim}_{A}^\theta F=0\ \mbox{if and only if}\ \dim_{A}^\theta F=0.$$
\end{theorem}
\noindent\textbf{Remark.}
The proof method of the equivalence that the upper spectrum is of full dimension is equivalent to the Assouad spectrum is of full dimension is from the method we used in the proof of Theorem 2.6(1) in \cite{WL}. In \cite{FY2}, Fraser and Yu proved that, in general, $ \overline{\dim}_{A}^\theta$ and $\dim_{A}^\theta$ are not necessarily equal.

\begin{example}\label{egb}
Fix any positive integer $k\geq 2$.
For any $n\in\mathbb{N}$, let $\delta_n:=\frac{1}{n^{k/(k-1)}}$ and $E_k:=\bigcup\limits_{n=1}^{+\infty}\{i\delta_n:i=0,1,\cdots, n\}.$
\noindent Let $E:=\bigcup\limits_{i=2}^{+\infty}(E_i+i)$,
where $A+b:=\{a+b:a\in A\}$ for any $A\subset\mathbb{R}$ and $b\in\mathbb{R}$.
\noindent It is easy to see that $E$ is a countable unbounded set.
We have
$$\dim_H E=0\ \mbox{and}\ \overline{\dim}_{GB} E=1.$$
\end{example}

\begin{proposition}\label{mgubde}
For any $F\subset\mathbb{R}^d$, we define
$$a(F):=\inf\bigg\{\sup\limits_{i\geq 1}\overline{\dim}_{GB} F_i: F\subset\bigcup\limits_{i=1}^{+\infty}F_i\bigg\}$$
and
$$b(F):=\inf\bigg\{\sup\limits_{i\geq 1}\overline{\dim}_{GB} F_i: F\subset\bigcup\limits_{i=1}^{+\infty}F_i,\ \mbox{where}\ F_i\ \mbox{is bounded for each}\ i\in \mathbb{N}\bigg\}.$$
Then
$$a(F)=b(F).$$
\end{proposition}

\begin{definition}
For any $F\subset\mathbb{R}^d$, the modified generalized upper box dimension of $F$, denoted by $\overline{\dim}_{MGB} F$, is defined by
$$\overline{\dim}_{MGB} F:=\inf\bigg\{\sup\limits_{i\geq 1}\overline{\dim}_{GB} F_i: F\subset\bigcup\limits_{i=1}^{+\infty}F_i\bigg\}.$$
\end{definition}

\begin{proposition} \label{mgubdpd}
For any $F\subset\mathbb{R}^d$, we have $$\dim_P F=\overline{\dim}_{MGB} F.$$
\end{proposition}

By Propositions \ref{gubdubd} and \ref{mgubdpd}, We know that the generalized upper box and upper box dimensions coincide for bounded sets. Furthermore, we also obtain that the modified generalized upper box dimension equals packing dimension. It is known that the modified upper box dimension equals packing dimension. So the generalized upper box dimension can be seen as a reasonable generalization of the upper box dimension.

\textbf{Organization.} In Sec \ref{gubd}, we will give the proof of results from Proposition \ref{asusl} to Example \ref{egb}. The proof of the remaining results can be found in Sec \ref{mgubd}.

\section{Generalized upper box dimension}\label{gubd}

\begin{proof}[Proof of Proposition \ref{asusl}]
By Proposition \ref{usasbp}(2), we have
$$\dim_{A}^\theta F \leq \overline{\dim}_{A}^\theta F$$
for each $\theta\in (0, 1)$, then
$$\limsup_{\theta\rightarrow 0}\dim_{A}^\theta F\leq\lim_{\theta\rightarrow 0}\overline{\dim}_{A}^\theta F.$$
For any $n\in \mathbb{N}$, there exists a $0<\theta_n<1$ such that
$$\lim_{n\rightarrow 0}\theta_n=0\ \mbox{and}\ |\overline{\dim}_{A}^{\theta_n} F-\lim_{\theta\rightarrow 0}\overline{\dim}_{A}^\theta F|<\frac{1}{2n}.$$
For $n$ above, by Proposition \ref{usas}, we have
$$\overline{\dim}_{A}^{\theta_n} F=\sup\limits_{0<\theta^{'}\leq \theta_n}\dim_{A}^{\theta^{'}} F,$$
then there exists a $0<\theta_n^{'}\leq \theta_n$ such that
$$|\overline{\dim}_{A}^{\theta_n} F-\dim_{A}^{\theta_n^{'}} F|<\frac{1}{2n}.$$
Then we have
$$\lim_{n\rightarrow 0}\theta_n^{'}=0$$
and
\begin{align*}
  |\dim_{A}^{\theta_n^{'}} F-\lim_{\theta\rightarrow 0}\overline{\dim}_{A}^\theta F| & \leq |\dim_{A}^{\theta_n^{'}} F-\overline{\dim}_{A}^{\theta_n} F|+|\overline{\dim}_{A}^{\theta_n} F-\lim_{\theta\rightarrow 0}\overline{\dim}_{A}^\theta F| \\
   & <\frac{1}{2n}+\frac{1}{2n} \\
   & =\frac{1}{n}.
\end{align*}
So
$$\lim_{n\rightarrow +\infty}\dim_{A}^{\theta_n^{'}} F=\lim_{\theta\rightarrow 0}\overline{\dim}_{A}^\theta F.$$
Furthermore, since
$$\lim_{n\rightarrow +\infty}\theta_n^{'}=0$$
and
$$\limsup_{\theta\rightarrow 0}\dim_{A}^\theta F\leq\lim_{\theta\rightarrow 0}\overline{\dim}_{A}^\theta F,$$
we have
$$\limsup_{\theta\rightarrow 0}\dim_{A}^\theta F=\lim_{\theta\rightarrow 0}\overline{\dim}_{A}^\theta F.$$

\end{proof}

\begin{proof}[Proof of Proposition \ref{gubdubd}]
By Proposition \ref{asub}, we have
$$\lim_{\theta\rightarrow 0}\dim_{A}^\theta F=\overline{\dim}_B F.$$
Thus,
$$\overline{\dim}_{GB} F=\limsup_{\theta\rightarrow 0}\dim_{A}^\theta F=\lim_{\theta\rightarrow 0}\dim_{A}^\theta F=\overline{\dim}_B F.$$
\end{proof}

\begin{proof}[Proof of Example \ref{gubr}]
Let $r_n=\delta_n$, $R_n=n\delta_n$ and $x_n=n$ for any positive integer $n\geq 2$.

\noindent For any $k\geq 2$, $0<\varepsilon<1$ and $C>0$, there exists a positive integer $N\geq k$ such that, for any $n>N$,
we have

$$R_n^k=(n\delta_n)^k=\left(\frac{n}{n^{1+1/(n-1)}}\right)^k=\left(\frac{1}{n^{1/(n-1)}}\right)^k
=\frac{1}{n^{k/(n-1)}}\geq \frac{1}{n^{n/(n-1)}}=\delta_n=r_n$$
and
$$\frac{1}{3}\left(\frac{R_n}{r_n}\right)^\varepsilon>C.$$
Then, for $0<r_n\leq R_n^k<R_n<1$ and $x_n\in E$, we have
\begin{align*}
  N_{r_n}(B(x_n,R_n)\cap E) & \geq \frac{n}{3} \\
   & =\frac{1}{3}\left(\frac{R_n}{r_n}\right) \\
   & =\frac{1}{3}\left(\frac{R_n}{r_n}\right)^\varepsilon\left(\frac{R_n}{r_n}\right)^{1-\varepsilon} \\
   & >C\left(\frac{R_n}{r_n}\right)^{1-\varepsilon}.
\end{align*}
So
$$\overline{\dim}_{A}^\frac{1}{k} E\geq 1-\varepsilon.$$
It follows from the arbitrariness of $\varepsilon$ and $k$ that
$$\overline{\dim}_{A}^\frac{1}{k} E=1\ \mbox{and}\ \overline{\dim}_{GB} E=1.$$
It is obvious that
$$\overline{\dim}_{GB}^* E=0.$$
\end{proof}

\begin{proof}[Proof of Theorem \ref{gubdasus}]
By Proposition \ref{ai}, for any $0<\theta_1<\theta_2<1$, we have
$$\dim_{A}^{\theta_1} F \leq  \left(\frac{1-\theta_2}{1-\theta_1}\right)\dim_{A}^{\theta_2} F+ \left(\frac{\theta_2-\theta_1}{1-\theta_1}\right)\dim_{A}^{\theta_1/\theta_2} F.$$
Then we have
$$\limsup_{\theta_1\rightarrow 0}\dim_{A}^{\theta_1} F \leq  (1-\theta_2)\dim_{A}^{\theta_2} F+ \theta_2\limsup_{\theta_1\rightarrow 0}\dim_{A}^{\theta_1/\theta_2} F,$$
then
$$\overline{\dim}_{GB} F \leq  (1-\theta_2)\dim_{A}^{\theta_2} F+ \theta_2\overline{\dim}_{GB} F,$$
so
$$\overline{\dim}_{GB} F\leq \dim_{A}^{\theta_2} F.$$
By Proposition \ref{ai}, for any $0<\theta_1<\theta_2<1$, we have
$$\left(\frac{1-\theta_2}{1-\theta_1}\right)\dim_{A}^{\theta_2} F \leq \dim_{A}^{\theta_1} F.$$
Then we have
$$(1-\theta_2)\dim_{A}^{\theta_2} F \leq \limsup_{\theta_1\rightarrow 0}\dim_{A}^{\theta_1} F.$$
So
$$\dim_{A}^{\theta_2} F\leq \frac{\overline{\dim}_{GB} F}{1-\theta_2}.$$
For any $0<\theta<1$, by Proposition \ref{usasbp}(2), we have
$$\dim_{A}^\theta F\leq \overline{\dim}_{A}^\theta F\leq \dim_{qA}F.$$
Thus, for $0<\theta<1$ and $F\subset\mathbb{R}^d$, we have
$$\overline{\dim}_{GB} F \leq \dim_{A}^\theta F\leq \min\left\{\frac{\overline{\dim}_{GB} F}{1-\theta}, \dim_{qA}F\right\}.$$
Since
$$\dim_{A}^\theta F\leq\frac{\overline{\dim}_{GB} F}{1-\theta}$$
for each $\theta\in (0, 1),$
we have
$$\dim_{A}^{\theta^{'}} F\leq\frac{\overline{\dim}_{GB} F}{1-\theta^{'}}\leq\frac{\overline{\dim}_{GB} F}{1-\theta}$$
for each $\theta^{'}\in (0, \theta].$

\noindent By Proposition \ref{usas}, we have
$$\overline{\dim}_{A}^\theta F=\sup\limits_{0<\theta^{'}\leq \theta}\dim_{A}^{\theta^{'}} F\leq \frac{\overline{\dim}_{GB} F}{1-\theta}.$$
By Proposition \ref{usasbp} and the definition of the generalized upper box dimension, we have
$$\overline{\dim}_{GB} F\leq \overline{\dim}_{A}^\theta F\leq \dim_{qA} F$$
for each $\theta\in (0, 1).$
So we have
$$\overline{\dim}_{GB} F \leq \overline{\dim}_{A}^\theta F\leq \min\left\{\frac{\overline{\dim}_{GB} F}{1-\theta}, \dim_{qA}F\right\}.$$

\end{proof}

\begin{proof}[Proof of Proposition \ref{gubdas}]
By Theorem \ref{gubdasus}, we have
$$\overline{\dim}_{GB} F \leq \dim_{A}^\theta F\leq \frac{\overline{\dim}_{GB} F}{1-\theta}$$
for each $\theta\in (0,1)$.
Thus we have
$$\overline{\dim}_{GB} F =\lim_{\theta\rightarrow 0}\dim_{A}^\theta F.$$
\end{proof}

\begin{proof} [Proof of Proposition \ref{gubdp}]
(1) By Proposition \ref{gubdas} and Proposition \ref{asbp}(1), we have
\begin{align*}
  \overline{\dim}_{GB} E & =\lim\limits_{\theta\rightarrow 0}\dim_{A}^\theta E \\
   & \leq \lim\limits_{\theta\rightarrow 0}\dim_{A}^\theta F \\
   & = \overline{\dim}_{GB} F.
\end{align*}
(2) By Proposition \ref{gubdas} and Proposition \ref{asbp}(2), we have
\begin{align*}
  \overline{\dim}_{GB} E\cup F & = \lim\limits_{\theta\rightarrow 0}\dim_{A}^\theta E\cup F\\
   & = \lim\limits_{\theta\rightarrow 0}\max\{\dim_{A}^\theta E, \dim_{A}^\theta F\}\\
   & = \lim\limits_{\theta\rightarrow 0}\frac{\dim_{A}^\theta E+\dim_{A}^\theta F+|\dim_{A}^\theta E-\dim_{A}^\theta F|}{2}\\
   & = \frac{\lim\limits_{\theta\rightarrow 0}\dim_{A}^\theta E+\lim\limits_{\theta\rightarrow 0}\dim_{A}^\theta F+|\lim\limits_{\theta\rightarrow 0}\dim_{A}^\theta E-\lim\limits_{\theta\rightarrow 0}\dim_{A}^\theta F|}{2} \\
   & = \frac{\overline{\dim}_{GB} E+\overline{\dim}_{GB} F+|\overline{\dim}_{GB} E-\overline{\dim}_{GB} F|}{2} \\
   & =\max\{\overline{\dim}_{GB} E, \overline{\dim}_{GB} F\}.
\end{align*}
(3) By Proposition \ref{gubdas} and Proposition \ref{asbp}(3), we have
\begin{align*}
  \overline{\dim}_{GB} E & =\lim\limits_{\theta\rightarrow 0}\dim_{A}^\theta E \\
   & = \lim\limits_{\theta\rightarrow 0}\dim_{A}^\theta F \\
   & = \overline{\dim}_{GB} F.
\end{align*}
(4) By Proposition \ref{gubdas} and Proposition \ref{asbp}(4), we have
\begin{align*}
  \overline{\dim}_{GB} F & =\lim\limits_{\theta\rightarrow 0}\dim_{A}^\theta F \\
   & = \lim\limits_{\theta\rightarrow 0}\dim_{A}^\theta \overline{F} \\
   & = \overline{\dim}_{GB} \overline{F}.
\end{align*}
(5) By Proposition \ref{gubdas} and Proposition \ref{asbp}(5), we have
\begin{align*}
  \overline{\dim}_{GB} E\times F & =\lim\limits_{\theta\rightarrow 0}\dim_{A}^\theta E\times F \\
   & \leq\lim\limits_{\theta\rightarrow 0}(\dim_{A}^\theta E+\dim_{A}^\theta F) \\
   & =\lim\limits_{\theta\rightarrow 0}\dim_{A}^\theta E+\lim\limits_{\theta\rightarrow 0}\dim_{A}^\theta F \\
   & =\overline{\dim}_{GB} E+\overline{\dim}_{GB} F.
\end{align*}
(6) By Proposition \ref{pdmubd}, Proposition \ref{asub}, Theorem \ref{gubdasus} and Proposition \ref{gubdp}(1), for any $0<\theta<1$, we have
\begin{align*}
  \dim_P F & =\inf\bigg\{\sup\limits_{i\geq 1}\overline{\dim}_B F_i: F=\bigcup\limits_{i=1}^{+\infty}F_i,\ \mbox{where}\ F_i\ \mbox{is bounded for each}\ i\in \mathbb{N}\bigg\} \\
   & \leq\inf\bigg\{\sup\limits_{i\geq 1}\dim_A^\theta F_i: F=\bigcup\limits_{i=1}^{+\infty}F_i,\ \mbox{where}\ F_i\ \mbox{is bounded for each}\ i\in \mathbb{N}\bigg\}  \\
   & \leq\inf\bigg\{\sup\limits_{i\geq 1}\frac{\overline{\dim}_{GB} F_i}{1-\theta}: F=\bigcup\limits_{i=1}^{+\infty}F_i,\ \mbox{where}\ F_i\ \mbox{is bounded for each}\ i\in \mathbb{N}\bigg\} \\
   & \leq \frac{\overline{\dim}_{GB} F}{1-\theta}.
\end{align*}
It follows from the arbitrariness of $\theta$ that
$$\dim_P F\leq \overline{\dim}_{GB} F.$$
By Theorem \ref{gubdasus} and Proposition \ref{usasbp}(2), we have $$\overline{\dim}_{GB} F\leq \dim_{qA} F\leq \dim_A F.$$

\end{proof}

\begin{proof}[Proof of Theorem \ref{gubdqad}]
By Theorem \ref{gubdasus}, we have $$\overline{\dim}_{GB} F\leq \dim_{qA} F.$$
Then if $$\dim_{qA} F=0,$$
we have
$$\overline{\dim}_{GB} F=0.$$
By Theorem \ref{gubdasus}, we have
$$\overline{\dim}_{GB} F \leq \overline{\dim}_{A}^\theta F\leq \frac{\overline{\dim}_{GB} F}{1-\theta}$$
for each $\theta\in (0, 1).$

\noindent Then if $$\overline{\dim}_{GB} F=0,$$ we have
$$\overline{\dim}_{A}^\theta F=0$$
for each $\theta\in (0, 1).$

\noindent It follows from the definition of the quasi-Assouad dimension that
$$\dim_{qA} F=\lim_{\theta\rightarrow 1}\overline{\dim}_{A}^\theta F=0.$$

\end{proof}

\begin{proof} [Proof of Theorem \ref{ua}]
For any $\theta\in (0,1)$, if $$\dim_{A}^\theta F=d,$$
then by Proposition \ref{usasbp}(2), we have
$$\overline{\dim}_{A}^\theta F\geq \dim_{A}^\theta F=d,$$
then
$$\overline{\dim}_{A}^\theta F=d.$$
Next we assume that
$$\overline{\dim}_{A}^\theta F=d.$$
For any $\theta\in (0,1)$, from Theorem \ref{gubdasus}, if $\overline{\dim}_{GB} F=d$, then we have
$$\dim_A^\theta F\geq \overline{\dim}_{GB} F=d,$$
then
$$\dim_A^\theta F=d.$$

\noindent If $\overline{\dim}_{GB} F<d$, then there exists $\theta_1\in (0, \theta)$ such that
$$\overline{\dim}_A^{\theta_1} F<d,$$
then it follows from Propositions \ref{usas} and \ref{usasbp} that
$$\overline{\dim}_A^\theta F=\sup\limits_{\theta'\in (0, \theta]}\dim_A^{\theta'} F=\sup\limits_{\theta'\in [\theta_1, \theta]}\dim_A^{\theta'} F=d.$$

\noindent By Proposition \ref{asc}, there exists $\theta_2\in [\theta_1, \theta]$ such that
$$\dim_A^{\theta_2} F=d=\dim_A F.$$
It follows from Proposition \ref{asa} that
$$\dim_A^\theta F=d.$$

For any $\theta\in (0,1)$, if $$\dim_{A}^\theta F=0,$$
by Theorem \ref{gubdasus}, we have
$$\overline{\dim}_{GB} F\leq \dim_{A}^\theta F=0,$$
so
$$\overline{\dim}_{GB} F=0.$$
By Theorem \ref{gubdasus}, we have
$$\overline{\dim}_{A}^\theta F\leq \frac{\overline{\dim}_{GB} F}{1-\theta},$$
so
$$\overline{\dim}_{A}^\theta F=0.$$
If
$$\overline{\dim}_{A}^\theta F=0,$$
by Proposition \ref{usasbp}(2), we have
$$\dim_{A}^\theta F\leq \overline{\dim}_{A}^\theta F=0,$$
then
$$\dim_{A}^\theta F=0.$$
\end{proof}

\begin{proof} [Proof of Example \ref{egb}]
Fix any positive integer $k\geq 2$.

\noindent For any positive integer $n\geq 2$, we have
$$n\delta_n<1$$
and
$$(n\delta_n)^k=\left(\frac{n}{n^{k/(k-1)}}\right)^k=\frac{1}{n^{k(k/(k-1)-1)}}=\frac{1}{n^{k/(k-1)}}=\delta_n.$$

For any $0<\varepsilon<1, C>0$, there exists a positive integer $N\geq 2$ such that, for any $n>N$,
we have
$$\frac{n^\varepsilon}{3}>C$$ and $$n\delta_n=\frac{n}{n^{k/(k-1)}}<1.$$
For any $n>N$, let $R_n=n\delta_n, r_n=\delta_n$ and $x_n=0.$

\noindent Then we have
$$0<r_n=R_n^k<R_n<1$$
and
$$x_n\in E_k,$$
then
\begin{align*}
  N_{r_n}(B(x_n,R_n)\cap E_k) & \geq \frac{n}{3}\\
   & =\frac{1}{3}\cdot\frac{R_n}{r_n} \\
   & = \frac{n^\varepsilon}{3}\left(\frac{R_n}{r_n}\right)^{1-\varepsilon}\\
   & >C\left(\frac{R_n}{r_n}\right)^{1-\varepsilon}.
\end{align*}
It follows from the arbitrariness of $\varepsilon$ that $\dim_{A}^{1/k}E_k=1.$

\noindent For any positive integer $i\geq 2,$ we have
$$1=\dim_{A}^{1/i}E_i=\dim_{A}^{1/i}(E_i+i).$$
It follows that
$$\overline{\dim}_{GB} E=\lim_{i\rightarrow +\infty}\dim_{A}^{1/i}E=1.$$
Since $E$ is countable, we have
$$\dim_H E=0.$$

\end{proof}

\section{Modified generalized upper box dimension} \label{mgubd}

\begin{proof} [Proof of Proposition \ref{mgubde}]
It is obvious that $$a(F)\leq b(F).$$
For any $(j_1, \cdots, j_d)\in \mathbb{Z}^d$, let $$Q_{j_1, \cdots, j_d}:=[j_1, j_1+1)\times\cdots\times [j_d, j_d+1).$$
For any sequence $\{F_i\}_{i=1}^{+\infty}$ of subsets of $\mathbb{R}^d$ with $F\subset\bigcup\limits_{i=1}^{+\infty}F_i$, we have
$$\sup\limits_{i\geq 1}\overline{\dim}_{GB} F_i\geq \sup\limits_{i\geq 1, (j_1, \cdots, j_d)\in \mathbb{Z}^d}\overline{\dim}_{GB} (F_i \cap Q_{j_1, \cdots, j_d}),$$
where $F\subset\bigcup\limits_{i\geq 1, (j_1, \cdots, j_d)\in \mathbb{Z}^d}(F_i \cap Q_{j_1, \cdots, j_d})$ and $F_i \cap Q_{j_1, \cdots, j_d}$ is bounded for each $i\in \mathbb{N}$ and $(j_1, \cdots, j_d)\in \mathbb{Z}^d$.

\noindent Then we have
\begin{align*}
  a(F) & =\inf\bigg\{\sup\limits_{i\geq 1}\overline{\dim}_{GB} F_i: F\subset\bigcup\limits_{i=1}^{+\infty}F_i\bigg\} \\
   & \geq \inf\bigg\{\sup\limits_{i\geq 1}\overline{\dim}_{GB} F_i: F\subset\bigcup\limits_{i=1}^{+\infty}F_i,\ \mbox{where}\ F_i\ \mbox{is bounded for each}\ i\in \mathbb{N}\bigg\} \\
   & =b(F).
\end{align*}

\end{proof}

\begin{proof} [Proof of Proposition \ref{mgubdpd}]
By Proposition \ref{pdmubd}, Proposition \ref{gubdubd} and Proposition \ref{mgubde}, we have
\begin{align*}
  \dim_P F & = \overline{\dim}_{MB} F\\
   & =\inf\bigg\{\sup\limits_{i\geq 1}\overline{\dim}_B F_i: F\subset\bigcup\limits_{i=1}^{+\infty}F_i,\ \mbox{where}\ F_i\ \mbox{is bounded for each}\ i\in \mathbb{N}\bigg\} \\
   & =\inf\bigg\{\sup\limits_{i\geq 1}\overline{\dim}_{GB} F_i: F\subset\bigcup\limits_{i=1}^{+\infty}F_i,\ \mbox{where}\ F_i\ \mbox{is bounded for each}\ i\in \mathbb{N}\bigg\} \\
   &  =\inf\bigg\{\sup\limits_{i\geq 1}\overline{\dim}_{GB} F_i: F\subset\bigcup\limits_{i=1}^{+\infty}F_i\bigg\} \\
   & =\overline{\dim}_{MGB} F.
\end{align*}

\end{proof}

\noindent\textbf{Acknowledgements}

L. Wang was supported by "Haizhou Bay Talents" Innovation Program KQ25011. W. Li was supported by NSFC No.12071148 and NSFC No.12471085.



\bibliographystyle{abbrvnat}

\bibliography{gubd}

\end{document}